\documentclass[11pt,a4paper]{article}
\usepackage{pb-diagram}
\usepackage{amscd,amsmath,amsthm}
\newtheorem{theorem}{Theorem}[section]
\newtheorem{lemma}[theorem]{Lemma}

\newtheorem{corollary}[theorem]{Corollary}
\newtheorem{proposition}[theorem]{Proposition}

\newtheorem{definition}{Definition}[section]

 \input amssym.def
\input amssym.tex

\def\bc{\begin{center}}
\def\ec{\end{center}}

\begin{document}
\title{The natural quiver of an artinian algebra\thanks{Project supported by the Program for New Century Excellent Talents
in University (No.04-0522) and the National Natural Science Foundation of China (No.10571153)}}
\author{Fang Li\thanks{fangli@cms.zju.edu.cn}$\;\;\;\;$ and$\;\;$ Lili Chen\thanks{lilychen0229@yahoo.com.cn} \\{Department of Mathematics, Zhejiang
University}\\{Hangzhou, Zhejiang 310027, China}} \maketitle
\begin{abstract}
The motivation of this paper is to study the natural quiver of an
artinian algebra, a new kind of quivers, as a tool independing
upon the associated basic algebra.

In \cite{Li}, the notion of the natural quiver of an artinian
algebra was introduced  and then was used to generalize the
Gabriel theorem for non-basic artinian algebras splitting over
radicals and non-basic finite dimensional algebras with
2-nilpotent radicals via pseudo path algebras and generalized path
algebras respectively.

In this paper, firstly we consider the relationship between the
natural quiver and the ordinary quiver of a finite dimensional
algebra. Secondly, the generalized Gabriel theorem is obtained
for radical-graded artinian algebras. Moreover, Gabriel-type
algebras are introduced to outline
 those artinian algebras satisfying the generalized Gabriel
 theorem here and in \cite{Li}. For such algebras, the uniqueness
 of the related generalized path algebra and quiver holds up to
 isomorphism in the case when the ideal is admissible.
  For an artinian algebra, there are two basic
algebras, the first is that associated to the algebra itself; the
second is that associated to the correspondent generalized path
algebra. In the final part, it is shown that for a Gabriel-type
artinian algebra, the first basic algebra is a quotient of the
second basic algebra.

In the end, we give an example of a skew group algebra in which
the relation between the natural quiver and the ordinary quiver is
discussed.

\end{abstract}


\section{Natural quiver and the relation with ordinary quiver}

Suppose that $A$ is a left artinian algebra over a field $k$, and
$r=r(A)$ is the radical of $A$. In this paper left artinian
algebras are written briefly as ``artinian algebras".

 Let $\{S_{1},S_{2},\cdots, S_{n}\}$ be the complete set of non-isomorphic simple $A$-modules of $A$.
 One can define
  a finite quiver $\Gamma_A$, called the {\em ordinary quiver} of $A$ as follows: $\Gamma_{0}=\{1,2,\cdots,
  n\}$, and the number $m_{ij}$ of arrows from $i$ to $j$ equals
  to the dimensional number dim$_{k}$Ext$_A$$(S_{i},S_{j})$.
  By \cite{A}, when $A$ is a finite-dimensional basic algebra over an
  algebraically closed field $k$ and $1_A=\varepsilon_{1}+\cdots
  +\varepsilon_{n}$ a decomposition of $1_A$ into a sum of primitive
  orthogonal idempotents.  Then, we can re-index $\{S_{1},S_{2},\cdots,
  S_{n}\}$ such that $S_{i}\cong A\varepsilon_{i}/r\varepsilon_{i}$, and moreover,
  dim$_{k}$Ext$_A$$(S_{i},S_{j})=$dim$_{k}$$(\varepsilon_{j}r/r^{2}\varepsilon_{i})$.
  Clearly, if $Q $ is a finite quiver without oriented cycles, the ordinary quiver of the path
  algebra $kQ$ is just $Q$.

Now, we introduce the so-called {\em natural quivers} from
artinian algebras.

 Write $A/r=\bigoplus_{i=1}^{s}\overline{A}_{i}$ where
$\overline{A}_{i}$ is a simple ideal of $A/r$ for each $i$. Then,
the algebra $r/r^2$ is an $A/r$-bimodule by
$\bar{a}\cdot(r/r^2)\cdot \bar{b}=arb/r^2$ for any $\bar
a=a+r,\bar b=b+r\in A/r$. Let $\;_{i}M_{j}=\overline{A}_{i}\cdot
r/r^{2}\cdot \overline{A}_{j}$, then $\;_{i}M_{j}$ is finitely
generated as $\overline{A}_{i}$-$\overline{A}_{j}$-bimodule for
each pair $(i,j)$.

For two artinian algebras $A$ and $B$,  the {\em rank} of a
finitely generated $A$-$B$-bimodule $M$ is defined as the least
cardinal number of the sets of generators. Clearly, for any
finitely generated $A$-$B$-bimodule, such rank always exists
uniquely.

Now we can associate with $A$ a quiver
$\Delta_A=(\Delta_{0},\Delta_{1})$, which is called the {\em
natural quiver} of $A$, in the following way. Let
$\Delta_{0}=\{1,\cdots,s\}$ as the set of vertices. For $i,j\in
\Delta_0$, let the number $t_{ij}$ of arrows from $i$ to $j$ in
$\Delta_A$ be the rank  of the finitely generated
$A_{j}$-$A_{i}$-bimodules $_{j}M_{i}$. Obviously, if
$_{j}M_{i}=0$, there are no arrows from $i$ to $j$.

The notion of natural quiver was firstly introduced in \cite{Li},
where the aim of the author is to use the generalized path algebra
from the natural quiver of an artinian algebra $A$ to characterize
$A$ through the generalized Gabriel theorem. In the further
research, one is motivated to study the representation of an
artinian algebra via the associated generalized path algebra or
pseudo path algebra but not the basic algebra of the artinian
algebra.

In order to clean the relation between the ordinary quiver and the
natural quiver of  an artinian algebra, it is necessary to note
that the natural quiver defined here is indeed  opposite to the
quiver defined in \cite{Li}. Now, we consider the relation between
the ordinary quiver and the natural quiver  of an artinian
$k$-algebra over an algebraically closed field $k$.

Clearly the number of the vertices in two quivers are equal since
$A$ and $A/r$ have the same simple modules, that is, we have $n=s$
 as above.

  When $A$ is a finite-dimensional basic algebra over an algebraically closed field $k$,  $A/r\cong\prod
k$ where the number of
  copies of $k$ equals the number of primitive
  orthogonal idempotents. As mentioned above,
  $dim_kExt^1_{A}(S_i,S_j)=dim_k(\varepsilon_jr/r^2\varepsilon_i)$
 where $S_i$, $S_j$ are the simple modules of $A$ corresponding to the primitive orthogonal idempotents $\varepsilon_i$,
 $\varepsilon_j$ respectively,
    which means that the number of arrows from $i$ to $j$ in the ordinary quiver $\Gamma_A$ of $A$
    is equal to that in the natural quiver $\Delta_A$ of $A$.
 Thus, we have:
 \begin{lemma}
For a finite-dimensional basic algebra $A$ over an algebraically
closed field $k$, the ordinary quiver and the natural quiver of
$A$ coincide. In particular, if $Q $ is a finite quiver without
oriented cycles, the ordinary quiver and the natural quiver of the
path algebra $kQ$ are both $Q$.
\end{lemma}

In order to discuss similarly for non-basic algebras, we introduce
the following notion:

Let $Q$ be a quiver and $Q'$ a sub-quiver of $Q$. If $(Q')_0=Q_0$
and for any vertices $i,j$,  there exist arrows from $i$ to $j$ in
$Q'$ if and only if there exist arrows from $i$ to $j$ in $Q$,
then we call this $Q'$ a {\em dense sub-quiver} of $Q$.

 When $A$ is over an
algebraically closed field $k$, by Proposition 7.4.4 in
\cite{Liu}, the relation:
$$t_{ij}\leq m_{ij}\leq n_in_jt_{ij}$$ holds where $n_i$ and $n_j$ are integers such that $\overline A_i\cong
M_{n_i}(k)$ and $\overline A_j\cong M_{n_j}(k)$. Trivially, if
each $n_i=1$, then $t_{ij}=m_{ij}$,  thus the ordinary quiver
$\Delta_A$ and $\Gamma_A$ the natural quiver of $A$ are coincided.
But, when some $n_i\not=1$, it is possible that $t_{ij}<m_{ij}$,
and $t_{ij}\not=0$ if and only if $m_{ij}\not=0$, which means
usually, $\Delta_A$ is a dense sub-quiver of $\Gamma_A$.

As  well-known, for an artinian algebra $A$, there is the
correspondent basic algebra $B$ and they are Morita-equivalent,
i.e. the module categories Mod$A$ and Mod$B$ are equivalent, which
follows that there is an equivalent functor $F$ such that
$Hom_A(S_i,S_j)\stackrel F{\cong}Hom_B(F(S_i),F(S_j))$
 for any simple modules $S_i$ and $S_j$ in Mod$A$. Moreover, $Ext^1_A(S_i,S_j)\cong Ext^1_B(F(S_i),F(S_j))$.
 It means the ordinary quiver of $A$ is the same with that of $B$.
 If $A$ is of finite dimension, its basic algebra is also of
 finite dimension. In the summary, we have:
  \begin{proposition}
Let $A$ be a finite dimensional algebra over a field $k$ with $r$
its radical and $A/r=\overline A_1\oplus\cdots\oplus\overline A_n
$ the direct sum of simple ideals,
 and  $B$ is the corresponding basic algebra of $A$.  Let $\Gamma_A$ and
$\Gamma_B$ be the  ordinary quivers of $A$ and $B$ respectively,
meanwhile $\Delta_A$ and $\Delta_B$  the natural quivers of $A$
and $B$ respectively. Then,

(i)$\;$ $\Gamma_A=\Gamma_B$;

(ii)$\;$  $\Gamma_B=\Delta_B$ if $k$ is algebraically closed;

(iii)$\;$  $\Delta_A$ is a dense sub-quiver of $\Gamma_A$, also of
$\Gamma_B$ and $\Delta_B$, if $k$ is algebraically closed.

  \end{proposition}

\section{Generalized Gabriel theorem in the radical-graded case}

The concepts of generalized path algebras were introduced  early
in \cite{CL} in order to find a generalization of path algebras so
as to obtain a generalized type of the Gabriel Theorem for
arbitrary finite dimensional algebras which would admit this
algebra to be isomorphic to a quotient algebra of such a
generalized path algebra. It is natural to ask how we look for a
generalized path algebra via the natural quiver to cover the
artinian algebra. Unfortunately, in general, as shown by the
counter-example in \cite{Li}, an artinian algebra with lifted
quotient may not be a
 homomorphic image of its correspondent $\mathcal A$-path-type tensor algebra.
 In this reason, the concepts of pseudo path algebras were introduced in \cite{Li} and it was shown that
  when the quotient algebra of an
artinian algebra can be lifted, the algebra is covered by a pseudo
path algebra via the natural quiver under an algebra homomorphism.

  However, there still exists  some special class of  artinian
 algebras which can be covered by their correspondent $\mathcal A$-path-type tensor
 algebras
 and equivalently by the generalized path algebras. This point can be
 seen in \cite{Li} from the generalized Gabriel theorem for a finite dimensional algebra with 2-nilpotent
 radical in the case it is splitting over its radical.

 In this section, we will give another class of
 artinian algebras which can be covered by the generalized path
 algebra via the natural quiver, that is, the generalized Gabriel
 theorem for this class of artinian algebras is true, too.  This
 class of artinian algebras are just the so-called radical-graded artinian
 algebra as follows.

For an artinian algebra $A$, let $r=$ rad$A$ be the radical of $A$
and the Loewy length $rl(A)=s$. Define gr$A=A/r\oplus
r/r^{2}\oplus \cdots \oplus r^{s-2}/r^{s-1}\oplus r^{s-1}$ as a
graded-algebra with multiplication
$(x+r^{i+1})(y+r^{j+1})=xy+r^{i+j+1}$ for $x\in r^i$, $y\in r^j$.
Trivially, this graded algebra is strict.

An artinian algebra $A$ is said to be {\em  radical-graded} if
$A=\bigoplus_{i\geq 0}A_{i}$ is strictly graded with $A_0$
semisimple. In this case, there is a minimal positive integer $t$
such that $A_i=0$ for all $i\geq t$ since $A$ is artinian. By this
definition, it is easy to see that for any artinian algebra $A$,
gr$A$ is always radical-graded. We have the following
characterization:

\begin{proposition} An artinian algebra $A$ is radical-graded if and only if $A\cong
$gr$A$. In this situation, $A_0\cong A/r$, $r\cong r/r^{2}\oplus
\cdots \oplus r^{s-2}/r^{s-1}\oplus r^{s-1}$ as algebras.
\end{proposition}

{\em Proof.}  ``$\Longleftarrow$" is trivial since gr$A$ is
radical-graded.

``$\Longrightarrow$": $\;$ Suppose that $A=\bigoplus_{i\geq
0}A_{i}$ is strictly graded with $A_{0}$ semisimple. Thus, there
is a minimal positive integer $t$ such that $A_{i}=0$ for $i\geq
t$. Write $r'=\bigoplus_{i\geq 1}A_{i}$, clearly it is an ideal of
$A$, $A/r'=A_{0}$ semisimple and $r'^{j}=\bigoplus _{i\geq
j}A_{i}$ which is zero when $j\geq t$ and hence $r'$ is nilpotent.
So $r'=r$ the radical of $A$, and clearly $r'^{i}/r'^{i+1}\cong
A_{i}$ for $i\geq 1$ which is zero when $i\geq t$, and $r'^l=0$ if
and only if $A_l=0$ for any $l$. Therefore $t=s=rl(A)$,
gr$A=A/r\oplus r/r^{2}\oplus \cdots \oplus r^{s-1}\cong A/r'\oplus
r'/r'^{2}\oplus \cdots \oplus r'^{s-1}=\bigoplus_{i\geq
0}A_{i}=A$. $\;\;\;\square$

From the above proposition, we get

\begin{corollary} For any artinian algebra $A$, gr(gr$A$)$\cong$gr$A$.
\end{corollary}

Now, we introduce briefly some notions about generalized path
algebras.

Let $Q=(Q_{0},Q_{1})$ be a quiver and $\mathcal{A}=\{A_{i}: i\in
Q_{0}\}$  a family of $k$-algebras $A_{i}$ with identity $e_{i}$,
indexed by the vertices of $Q$. The elements $a_{j}\neq 0$ of
$\bigcup_{i\in Q_{0}}A_{i}$ are called  {\em $\mathcal{A}$-paths
of length zero}, with starting vertex $s(a_{j})$ and the ending
vertex $e(a_{j})$ are both $j$. For each $n\geq 1$, an {\em
$\mathcal{A}$-path $P$ of length $n$} is given by
$a_{1}\beta_{1}a_{2}\beta_{2}\cdot\cdot\cdot
a_{n}\beta_{n}a_{n+1}$, where
$(s(\beta_{1})|\beta_{1}\cdot\cdot\cdot\beta_{n}|e(\beta_{n}))$ is
a path in $Q$ of length $n$, for each $i=1,...,n$, $0\neq a_{i}\in
A_{s(\beta_{i})}$ and $0\neq a_{n+1}\in A_{e(\beta_{n})}$.
  $s(\beta_{1})$ and $e(\beta_{n})$ are also called respectively  the starting vertex and the ending vertex of
 $P$. Write $s(P)=s(\alpha_{1})$ and $e(P)=e(\alpha_{n})$.
 Now, consider the quotient $R$ of
the $k$-linear space with basis the set of all $\mathcal{A}$-paths
 by the subspace generated by all the elements of the form
\[
a_{1}\beta_{1}\cdot\cdot\cdot\beta_{j-1}(a_{j}^{1}+\cdot\cdot\cdot+a_{j}^{m})\beta_{j}a_{j+1}\cdot\cdot\cdot\beta_{n}a_{n+1}
-\sum_{l=1}^{m}a_{1}\beta_{1}\cdot\cdot\cdot\beta_{j-1}a_{j}^{l}\beta_{j}a_{j+1}\cdot\cdot\cdot\beta_{n}a_{n+1}
\]
where
$(s(\beta_{1})|\beta_{1}\cdot\cdot\cdot\beta_{n}|e(\beta_{n}))$ is
a path in $Q$ of length $n$, for each $i=1,...,n$, $a_{i}\in
A_{s(\beta_{i})}$, $a_{n+1}\in A_{e(\beta_{n})}$ and $a_{j}^{l}\in
A_{s(\beta_{j})}$ for $l=1,...,m$. In $R$,  given two elements
 $[a_{1}\beta_{1}a_{2}\beta_{2}\cdot\cdot\cdot a_{n}\beta_{n}a_{n+1}]$ and
 $[b_{1}\gamma_{1}b_{2}\gamma_{2}\cdot\cdot\cdot b_{n}\gamma_{n}b_{n+1}]$, define the multiplication as follows:
\\
\\
$ [a_{1}\beta_{1}a_{2}\beta_{2}\cdot\cdot\cdot
a_{n}\beta_{n}a_{n+1}]\cdot[b_{1}\gamma_{1}b_{2}\gamma_{2}\cdot\cdot\cdot
b_{n}\gamma_{n}b_{n+1}]\\
\\
=\left\{\begin{array}{ll}
 [a_{1}\beta_{1}a_{2}\beta_{2}\cdot\cdot\cdot
a_{n}\beta_{n}(a_{n+1}b_{1})\gamma_{1}b_{2}\gamma_{2}\cdot\cdot\cdot
b_{n}\gamma_{n}b_{n+1}],
  &  \mbox{if $a_{n+1}, b_{1}\in A_{i}$ for the same $i$}\\
0,  &  \mbox{otherwise}
\end{array}
\right. $ \\
It is easy to check that the above multiplication is well-defined
 and makes $R$ to become a $k$-algebra. This algebra $R$
defined above is called an {\em $\mathcal{A}$-path algebra} of $Q$
respecting to $\mathcal{A}$, or generally  {\em generalized path
algebras}. Denote it by $R=k(Q,\mathcal{A})$. Clearly, $R$ is an
$A$-bimodule, where $A=\oplus_{i\in Q_0}A_{i}$.

A generalized path algebra $k(Q,\mathcal{A})$ is said to be {\em
normal} if all algebras $A_i$ ($i\in Q_0$) are simple algebras for
$\mathcal A=\{A_i:\; i\in Q_0\}$.

Associated with the pair $(A,\;_{A}M_{A})$ for a $k$-algebra $A$
and an $A$-bimodule $M$, we write the $n$-fold $A$-tensor product
$M\otimes_{A}M\otimes\cdot\cdot\cdot\otimes_{A}M$ as $M^{n}$.
Writing $M^{0}=A$, then $T(A,M)=A\oplus M\oplus
M^{2}\oplus\cdot\cdot\cdot\oplus M^{n}\oplus\cdot\cdot\cdot $
becomes a $k$-algebra with multiplication induced by the natural
$A$-bilinear maps $M^{i}\times M^{j}\rightarrow M^{i+j}$ for
$i\geq 0$ and $j\geq 0$.  $T(A,M)$ is called the {\em tensor
algebra} of $M$ over $A$.

Define a special class of tensor algebras so as to characterize
generalized path algebras. An {\em $\mathcal{A}$-path-type tensor
algebra} is defined to be the tensor algebra $T(A,M)$ satisfying
that (i) $A=\bigoplus_{i\in I}A_{i}$ for a family of $k$-algebras
$\mathcal{A}=\{A_{i}: i\in I\}$, (ii) $M=\bigoplus_{i,j\in I}$$
_{i}M_{j}$ where $_{i}M_{j}$ are finitely generated
$A_{i}$-$A_{j}$-bimodules
 for all $i$ and $j$ in $I$ and $A_{k}\cdot_{i}M_{j}=0$
if $k\not=i$ and
 $_{i}M_{j}\cdot A_{k}=0$ if $k\not=j$.
 A {\em free
$\mathcal{A}$-path-type tensor algebra} is the
$\mathcal{A}$-path-type tensor algebra $T(A,M)$ whose each
finitely generated $A_{i}$-$A_{j}$-bimodule $_{i}M_{j}$ for $i$
and $j$ in $I$ is a free bimodule with a basis and the cardinality
of this basis is equal to the rank of $_{i}M_{j}$ as a finitely
generated $A_{i}$-$A_{j}$-bimodule.

In an $\mathcal A$-path algebra $k(Q,\mathcal{A})$, let
$A=\bigoplus_{i\in Q_{0}}A_{i}$. For any $i$, $j$, let
$_{i}M^{F}_{j}$ be the free $A_{i}$-$A_{j}$-bimodule with basis
given by the arrows from $j$ to $i$. Then the number of free
generators in the basis is the rank of $_{i}M^{F}_{j}$ as a
finitely generated bimodule. Define $A_{k}\cdot_{i}M^{F}_{j}=0$ if
$k\not=i$ and
 $_{i}M^{F}_{j}\cdot A_{k}=0$ if $k\not=j$. Then $M^{F}=\bigoplus_{j\rightarrow i}$$ _{i}M^{F}_{j}$ is
  an $A$-bimodule. We  get the unique free $\mathcal{A}$-path-type tensor algebras $T(A,M^{F})$.

Conversely, given an $\mathcal{A}$-path-type tensor algebra
$T(A,M)$ with  $\mathcal{A}=\{A_{i}: i\in I\}$ and finitely
generated $A_{i}$-$A_{j}$-bimodules $_{i}M_{j}$ for $i,j\in I$
such that $A=\bigoplus_{i\in I}A_{i}$, $M=\bigoplus_{i,j\in I}$$
_{i}M_{j}$, $A_{k}\cdot_{i}M_{j}=0$ if $k\not=i$ and
 $_{i}M_{j}\cdot A_{k}=0$ if $k\not=j$. Trivially, $_{i}M_{j}=A_{i}MA_{j}$.
 Let $r_{ij}$ be the rank of $_{j}M_{i}$.
One can associate with $T(A,M)$ a quiver $Q=(Q_{0},Q_{1})$, called
{\em the quiver of $T(A,M)$}, via
 $Q_{0}=I$ as the set of vertices and for $i,j\in I$, $r_{ij}$ as the number
of arrows from $i$ to $j$ in $Q$. Its $\mathcal{A}$-path algebra
$k(Q,\mathcal{A})$ is called {\em the corresponding
$\mathcal{A}$-path algebra of $T(A,M)$}. By definition, the quiver
of $T(A/r,r/r^{2})$ is just $\Delta_A$.

From the above discussion, every $\mathcal{A}$-path-type tensor
algebra $T(A,M)$ can be used to construct its corresponding
$\mathcal{A}$-path algebra $k(Q,\mathcal{A})$; but, from this
$\mathcal{A}$-path algebra $k(Q,\mathcal{A})$, we can get uniquely
the free $\mathcal{A}$-path-type tensor algebra $T(A,M^{F})$. In
summary, we have the following in \cite{Li}:

 \begin{lemma}
(i)\; Every $\mathcal{A}$-path-type tensor algebra $T(A,M)$ can be
used to construct uniquely the free $\mathcal{A}$-path-type tensor
algebra $T(A,M^{F})$. There is a surjective $k$-algebra morphism
$\pi$: $T(A,M^{F})\rightarrow T(A,M)$ such that
$\pi(_{i}M^{F}_{j})=\;$$_{i}M_{j}$ for any $i,j\in I$;

(ii)\; Let $T(A,M^{F})$ be the free $\mathcal{A}$-path-type tensor
algebra built by a $\mathcal{A}$-path algebra $k(Q,\mathcal{A})$.
Then there is a $k$-algebra isomorphism $\widetilde{\phi}$:
$T(A,M^{F})\rightarrow k(Q,\mathcal{A})$
 such that for any $t\geq 1$,
$\widetilde{\phi}(\bigoplus_{j\geq
 t}(M^{F})^{j})= J^{t}$;

(iii)\; Let $T(A,M)$ be an $\mathcal{A}$-path-type tensor algebra
with the corresponding  $\mathcal{A}$-path algebra
$k(Q,\mathcal{A})$. Then there is a surjective $k$-algebra
homomorphism $\widetilde{\varphi}$: $k(Q,\mathcal{A})\rightarrow
T(A,M)$ such that for any $t\geq
1$,$\widetilde{\varphi}(J^{t})=\bigoplus_{j\geq t}M^j$.

 Here $ J$ denotes the ideal generated by  all
 $\mathcal{A}$-paths of length $1$ in $k(Q,\mathcal{A})$.
\end{lemma}

In the sequel, we always denote by $J$ the ideal generated by all
generalized paths of length one in the discussed generalized path
algebras. When $Q$ is admissible, i.e. is acyclic, $J$ is just the
radical of a normal generalized path algebra $k(Q,\mathcal{A})$
(see \cite{CL}).

 A {\em
 relation} $\sigma$ on an  $\mathcal{A}$-path
algebra $k(\Delta,\mathcal{A})$ is a $k$-linear combination of
some $\mathcal{A}$-paths $P_{i}$ with the same starting vertex and
the same ending vertex, that is,
$\sigma=k_{1}P_{1}+\cdot\cdot\cdot+k_{n}P_{n}$ with $k_{i}\in k$
and $s(P_{1})=\cdot\cdot\cdot=s(P_{n})$ and
$e(P_{1})=\cdot\cdot\cdot=e(P_{n})$.  If
$\rho=\{\sigma_{t}\}_{t\in T}$ is a set of relations on
 $k(\Delta,\mathcal{A})$, the pair
 $(k(\Delta,\mathcal{A}),\rho)$ is called an
{\em $\mathcal{A}$-path algebra with relations}. Associated with
 $(k(\Delta,\mathcal{A}),\rho)$ is the quotient $k$-algebra
 $k(\Delta,\mathcal{A},\rho)\stackrel{\rm
def}{=}k(\Delta,\mathcal{A})/\langle\rho\rangle$, where
$\langle\rho\rangle$ denotes the ideal in $k(\Delta,\mathcal{A})$
generated by the set of relations $\rho$. When the length
$l(P_{i})$ of each $P_{i}$ is at least $j$, it holds
 $\langle\rho\rangle\subset J^{j}$.

Now, let $M=r/r^2$ as $A/r$-bimodule,
$\;_{i}M_{j}=\overline{A}_{i}\cdot r/r^{2}\cdot \overline{A}_{j}$,
then $\;_{i}M_{j}$ is finitely generated as
$\overline{A}_{i}$-$\overline{A}_{j}$-bimodule for each pair
$(i,j)$ and $M=\bigoplus_{i,j}\;_iM_{j}$.
 Thus, we get the tensor algebra
$T(A/r,r/r^2)$ and the corresponding generalized path algebra
$k(\Delta_{A},\mathcal{A})$ from the natural quiver $\Delta_A$ of
$A$.

A set of some $\mathcal A$-paths or their linear combinations in
$k(Q,\mathcal A)$ is said to be  {\em $\mathcal A$-finite} if all
$\mathcal A$-paths in this set are constructed from a finite
number of paths in $Q$ with elements of $\bigcup_{i\in Q_0}A_i$. A
quotient or an ideal of $k(Q,\mathcal A)$ is said to be {\em
$\mathcal A$-finitely generated} if it is generated by an
$\mathcal A$-finite set.

The following is the main result in this section:
\begin{theorem} {\em (Generalized Gabriel Theorem in radical-graded case)}$\;\;$Assume that $A$ is a radical-graded artinian
$k$-algebra. Then, there is an $\mathcal A$-finite set $\rho$ of
relations of $k(\Delta_{A},\mathcal{A})$ such that
 $A\cong k(\Delta_{A},\mathcal{A})/\langle\rho\rangle$ with
 $ J^{s}\subset\langle\rho\rangle\subset J$
 for some positive integer $s$.

 \end{theorem}
{\em Proof}:$\;$ Let $r$ be the radical of $A$ with the Loewy
length $rl(A)=s+1$. Since $A$ is radical-graded, we have $A\cong
A/r\oplus r/r^2\oplus r^2/r^3\oplus\cdots\oplus r^{s-1}/r^s\oplus
r^s$.
 Thus, $r\cong r/r^2\oplus r^2/r^3\oplus\cdots\oplus r^{s-1}/r^s\oplus r^s$ and $A\cong A/r\oplus r$ as algebras.

 Write $A/r=\oplus_{i=1}^s\overline{A}_i$ with simple ideals $\overline A_i$ for all $i$. Then, we have the $\mathcal A$-path type tensor algebra $T(A/r, r/r^2)$ with
 $\mathcal A=\{\overline A_i: i=1\cdots s\}$. Firstly, we can find a surjective morphism of algebras from $T(A/r, r/r^2)$ to $A$. In fact, for $r^m/r^{m+1}$,
  define $f_m:\;r/r^2\otimes_{A/r}\cdots\otimes_{A/r}r/r^2$ (with $m$ copies of $r/r^2$) $\longrightarrow r^m/r^{m+1}$ satisfying that
  $f_m(\overline{x}_1\otimes\cdots\otimes\overline{x}_m)=\overline{x_1\cdots x_m}$ where $\overline{x}_i\in r/r^2$ for $x_i\in r$ and
 $\overline{x_1\cdots x_m}\in r^m/r^{m+1}$. It is easy to see that $f_m$ is well-defined as a morphism of $A/r$-$A/r$-bimodules and trivially,
 $f_m$ is surjective. Then,  $f=id_{A/r}\oplus f_1\oplus\cdots\oplus f_m\oplus\cdots$ is a surjective algebra morphism from $T(A/r, r/r^2)$ to $A$,
   where $f_m=0$ when $m\geq s+1$.

 Moreover, by Lemma 2.3,
 there is a surjective $k$-algebra homomorphism
$\widetilde{\varphi}$: $k(\Delta_A,\mathcal{A})\rightarrow
T(A/r,r/r^{2})$
 such that for any $t\geq 1$,$\widetilde{\varphi}(J^{t})=\bigoplus_{j\geq t}(r/r^2)^{\otimes j}$, where $(r/r^2)^{\otimes j}$
 denotes $r/r^{2}\otimes_{A/r}r/r^{2}\otimes_{A/r}\cdot\cdot\cdot\otimes_{A/r}r/r^{2}$
with $j$ copies of $r/r^{2}$.
 Then, $f\widetilde{\varphi}$: $k(\Delta_A,\mathcal{A})\rightarrow A$ is a surjective algebra morphism. Therefore, for the kernel $I=ker(f\widetilde{\varphi})$, we obtain
 $k(\Delta_A,\mathcal{A})/I\cong A$.

Now, we prove that $\oplus_{j\geq rl(A)}(r/r^{2})^{\otimes
j}\subset Kerf\subset\oplus_{j\geq 2}(r/r^{2})^{\otimes j}$. In
fact, by the definition, $f_1=id_{r/r^2}$, so $f|_{A/r\oplus
r/r^2}=id_{A/r}\oplus f_1:\;A/r\oplus r/r^2\longrightarrow A$ is a
monomorphism with image intersecting $r^2$ trivially. It follows
that $Kerf\subset\oplus_{j\geq 2}(r/r^{2})^{\otimes j}$. On the
other hand, $f((r/r^{2})^{\otimes j})=0$ for $j\geq rl(A)$ since
$r^{j}=0$ in this case. Therefore we get $\oplus_{j\geq
rl(A)}(r/r^{2})^{\otimes j}\subset Kerf$.

 But, by Lemma 2.3, for $t=rl(A)$ and $t=2$ respectively,  $\widetilde{\varphi}(J^{rl(A)})=\oplus_{j\geq rl(A)}(r/r^{2})^{\otimes j}$ and
$\widetilde\varphi( J^{2})=\oplus_{j\geq 2}(r/r^{2})^{\otimes j}$.
So, $\widetilde\varphi(J^{rl(A)})\subset
Kerf\subset\widetilde\varphi( J^{2})$.

Then, we prove $
J^{t}\subset\widetilde\varphi^{-1}\widetilde\varphi(J^{t})\subset
J^{t}+\widetilde\phi(\oplus_{j\leq t-1}((r/r^{2})^F)^{\otimes
j})\cap\widetilde\phi(Ker\pi)$ for $t\geq 1$, where
$\widetilde\varphi$ is that in Lemma 2.3. Trivially, $
J^{t}\subset\widetilde\varphi^{-1}\widetilde\varphi(J^{t})$. On
the other hand, $\widetilde\varphi=\pi\widetilde\phi^{-1}$ and
then $\widetilde\varphi^{-1}=\widetilde\phi\pi^{-1}$. By Lemma
2.3(iii), $\widetilde\varphi( J^{t})=\oplus_{j\geq
t}(r/r^{2})^{\otimes j}$. From the definition of $\pi$ in Lemma
2.3, it can be seen that $\pi^{-1}(\oplus_{j\geq
t}(r/r^{2})^{\otimes j})\subset\oplus_{j\geq
t}((r/r^{2})^{F})^{\otimes j}+(\oplus_{j\leq
t-1}((r/r^{2})^{F})^{\otimes j})\cap$Ker$\pi$. Thus, by Lemma 2.3,
we have

 $\widetilde\varphi^{-1}\widetilde\varphi(J^{t})=\widetilde\phi\pi^{-1}(\oplus_{j\geq
t}(r/r^{2})^{j})\subset\widetilde\phi(\oplus_{j\geq
t}((r/r^{2})^F)^{\otimes j})+\widetilde\phi(\oplus_{j\leq
t-1}((r/r^{2})^F)^{\otimes j})\cap\widetilde\phi(Ker\pi)\\
= J^{t}+\widetilde\phi(\oplus_{j\leq t-1}((r/r^{2})^F)^{\otimes j})\cap\widetilde\phi(Ker\pi)$. \\
  Hence,

$
J^{rl(A)}\subset\widetilde\varphi^{-1}\widetilde\varphi(J^{rl(A)})
\subset\widetilde\varphi^{-1}(Ker{f})=I\subset\widetilde\varphi^{-1}\widetilde\varphi( J^{2})\\
\subset
J^{2}+\widetilde\phi(\oplus_{j\leq 1}((r/r^{2})^F)^{\otimes j})\cap\widetilde\phi(Ker\pi) \subset J^{2}+ J\cap\widetilde\phi(Ker\pi)$.\\
But,
$\widetilde\phi(Ker\pi)=\widetilde\phi(\pi^{-1}(0))=\widetilde\varphi^{-1}(0)=Ker\widetilde\varphi$.
Then, $ J^{rl(A)}\subset I\subset J^{2}+ J\cap
Ker\widetilde\varphi\subset J.$

Lastly, we present $I$ through an $\mathcal A$-finite set of
relations on $k(\Delta_A,\mathcal{A})$. $ J^{rl(A)}$ is the ideal
$\mathcal A$-finitely generated in $k(\Delta_A,\mathcal{A})$ by
all $\mathcal{A}$-paths of length $rl(A)$.
$k(\Delta_A,\mathcal{A})/J^{rl(A)}$ is generated $\mathcal
A$-finitely, under the meaning of isomorphism, by all
$\mathcal{A}$-paths of length less than $rl(A)$, so as well as
$I/J^{rl(A)}$ as a $k$-subspace. Then, $I$ is an $\mathcal
A$-finitely generated ideal in $k(\Delta_A,\mathcal{A})$. Assume
$\{\sigma_{l}\}_{l\in\Lambda }$ is a set of $\mathcal A$-finite
generators for the ideal $I$. For the identity $\overline{1}$ of
$A/r$, we have the decomposition of orthogonal idempotents
$\overline{1}=\overline{e}_{1}+\cdot\cdot\cdot+\overline{e}_{s}$,
where $\overline{e}_{i}$ is the identity of $\overline{A}_{i}$.
Then $\sigma_{l}=\overline{1}\sigma_{l}\overline{1}=\sum_{1\leq
i,j\leq s}\overline{e}_{i}\sigma_{l}\overline{e}_{j}$. Obviously,
 $\overline{e}_{i}\sigma_{l}\overline{e}_{j}$ can be expanded as a
 $k$-linear combination of some such $\mathcal{A}$-paths which
 have the same starting vertex $j$ and the same ending vertex $i$. So, $\sigma^{ilj}=\overline{e}_{i}\sigma_{l}\overline{e}_{j}$
 is a relation on the $\mathcal{A}$-path algebra
 $k(\Delta_A,\mathcal{A})$. Moreover, $I$ is generated by all $\sigma^{ilj}$
 due to $\sigma_{l}=\sum_{i,j}\sigma^{ilj}$. Therefore, for
 $\rho=\{\sigma^{ilj}: 1\leq i,j\leq s,\; l\in\Lambda\}$, we get
 $I=\langle\rho\rangle$.
  Hence
  $k(\Delta_A,\mathcal{A},\rho)=k(\Delta_A,\mathcal{A})/\langle\rho\rangle\cong A$
  with
  $\langle\rho\rangle=Ker(f\widetilde\varphi)$ and $J^{rl(A)}\subset\langle\rho\rangle
  \subset J^{2}+ J\cap Ker\widetilde\varphi\subset J$.
  $\;\;\;\square$

 The uniqueness of the correspondent generalized path algebra and natural quiver of a radical-graded artinian
 algebra holds
  up to isomorphism if the ideal $\langle\rho\rangle$ is restricted into $J^2$. That is, if there exists another quiver and its
 related generalized path algebra such that the same isomorphism relation as in Theorem 2.4 is satisfied, then this quiver and
 related
 generalized path algebra are just respectively the natural quiver and the corresponding one of
 the radical-graded artinian algebra. This can be seen as a
 special case of the uniqueness of the so-called Gabriel-type algebras, see Theorem 3.3 in the next section.

\section{Two basic algebras from an artinian algebra }

For an artinian algebra $A$, write $A/r=\bigoplus_{i=1}^{s}
\overline{A_{i}}$ with simple ideals $\overline{A_{i}}$, we get
$k(\Delta_{A}, \mathcal{A})$ where $\Delta_{A}$ is the natural
quiver of $A$ and
$\mathcal{A}=\{\overline{A_{i}}:i=1,2,\cdots,s\}$.

It is known that the associated basic algebra $B$ which is
Morita-equivalent to $A$ is important for representations of $A$.
In order to realize our approach, it is valid to consider the
associated basic algebra $C$ of the generalized path algebra
$k(\Delta_A,\mathcal A)$ of the natural quiver $\Delta_A$ of $A$
and moreover, the relationship between $B$ and $C$.

However, in general, the generalized path algebra
$k(\Delta_A,\mathcal A)$ is not an artinian algebra, e.g. when the
natural quiver $\Delta_A$ contains an oriented cycle. So,
$k(\Delta_A,\mathcal A)$ has not the so-called related basic
algebra under the meaning of ``artinian" such that they are
Morita-equivalent each other. In this reason,
 $C$ is different from that for artinian
algebras.

A {\em complete set of non-isomorphic primitive orthogonal
idempotents} of $A$ is a set of primitive orthogonal idempotents
$\{\varepsilon_i:i\in I\subset(\Delta_A)_0\}$ such that
$A\varepsilon_{i}\not\cong A\varepsilon_{j}$ as left $A$-modules
for any $i\not=j$ in $I$ and for each primitive idempotent
$\varepsilon_s$ the module $A\varepsilon_s$ is isomorphic to one
of the modules $A\varepsilon_i$ ($i\in I$).

Every indecomposable projective module $P$ is decided by a
primitive idempotent $e_i$, that is, $P\cong A\varepsilon_i$ for
some $i$. And, there exists a bijective correspondence between the
iso-classes of indecomposable projective modules and the
iso-classes of simple modules. The set of the latter is equal to
the vertex set $(\Gamma_{A})_{0}$ of the ordinary quiver
$\Gamma_A$ of $A$, and then to the vertex set $(\Delta_{A})_{0}$
of the natural quiver $\Delta_A$ of $A$. Hence,
$I=(\Delta_{A})_{0}$. Let each $P_i$ be chosen as a representative
from the iso-class of indecomposable projective module
$A\varepsilon_i$ and let $i$ run over the vertex set
$(\Delta_A)_0$. Then the basic algebra $B$ of $A$ is given by
$B=End(\coprod_{i\in(\Delta_A)_0}P_i)\cong\bigoplus_{i,j\in(\Delta_A)_0}Hom_A(P_i,P_j)\cong\bigoplus_{i,j\in(\Delta_A)_0}\varepsilon_iA\varepsilon_j$.

\begin{lemma} Let $A$ be an artinian algebra. Then the complete set of non-isomorphic primitive orthogonal idempotents of $A$,
$A/r$ ($r$ is the radical of $A$) and $k(\Delta_A,\mathcal A)$ are
the same, whose cardinality is equal to that of the vertex set of
the natural quiver of $A$.

\end{lemma}

{\em Proof}: Let $\overline \varepsilon_i$ be the image of
$\varepsilon_i$ under the canonical homomorphism from $A$ to
$A/r$. Since $A$ and $A/r$ have the same simple modules,
$\{\overline{\varepsilon_i}:i\in (\Delta_A)_0\}$ is a complete set
of non-isomorphic primitive orthogonal idempotents of $A/r$. But
the idempotents of $k(\Delta_A,\mathcal A)$ must have length zero,
hence $\{\overline{\varepsilon_i}:i\in (\Delta_A)_0\}$ is also a
complete set of non-isomorphic primitive orthogonal idempotents of
$k(\Delta_A,\mathcal A)$. $\;\;\;\square$

As discussed before Lemma 3.1,  the basic algebra $C$ satisfies
$$C=End(\coprod_{i\in(\Delta_A)_0}k(\Delta_A,\mathcal
A)\overline\varepsilon_i)\cong\bigoplus_{i,\,j\in(\Delta_A)_0}\overline
\varepsilon_ik(\Delta_A,\mathcal A)\overline\varepsilon_j.$$ Then,
we get the following:

\begin{proposition}
For an artinian algebra $A$ over a field $k$ with the natural
quiver $\Delta_A$, let $\{\varepsilon_i:i\in(\Delta_A)_0\}$ be the
complete set of non-isomorphic primitive orthogonal idempotents of
$A$. Denote by $\overline \varepsilon_i$ the image of
$\varepsilon_i$ under the canonical morphism from $A$ to $A/r$.
Then,

(i) the basic algebra $B$ of $A$ is isomorphic to
$\bigoplus_{i\in(\Delta_A)_0}\varepsilon_iA\varepsilon_j$;

(ii) the basic algebra $C$ of the associated generalized path
algebra $k(\Delta_A,\mathcal A)$ of $A$ is isomorphic  to
$\bigoplus_{i\in(\Delta_A)_0}\overline
\varepsilon_ik(\Delta_A,\mathcal A)\overline \varepsilon_j$.
\end{proposition}

As we have said, $k(\Delta_A,\mathcal A)$ is not artinian when
$\Delta_A$ has an oriented cycle. Hence, we cannot affirm whether
$C$ is Morita equivalent to $k(\Delta_A,\mathcal A)$ in general.
But, $C$ is still decided uniquely by $k(\Delta_A,\mathcal A)$
 and then by $A$.

For an arbitrary artinian algebra $A$, we still cannot  obtain
 the explicit relation between two basic algebras $B$ and $C$
depending upon Proposition 3.2.
 However, for the following special case, that is, for the so-called {\em Gabriel-type algebras},
 we will give an exact conclusion for the two basic algebras.

\begin{definition}
Let $A$ be an artinian algebra over a field $k$ and
$k(\Delta_A,\mathcal A)$ its associated normal generalized path
algebra. If there exists an ideal $I$ of $k(\Delta_A,\mathcal A)$
such that $A\cong k(\Delta_A,\mathcal A)/I$, then we say $A$ to be
of {\em Gabriel-type}.
\end{definition}

Since in \cite{Li}, we have the Generalized Gabriel Theorem for a
finite dimensional algebra $A$ with 2-nilpotent radical $r=r(A)$
in the case $A$ is splitting over $r$, that is,
  $A\cong k(\Delta_{A},\mathcal{A})/\langle\rho\rangle$ with
$ J^{2}\subset\langle\rho\rangle\subset J^{2}+ J\cap$
\textrm{Ker}$\widetilde{\varphi}$ where $\langle\rho\rangle$ is an
ideal generated by the set of relations $\rho$ of
$k(\Delta,\mathcal{A})$ and $\widetilde\varphi$ is that in Lemma
2.3. It means that any such finite dimensional algebra is always
of Gabriel-Type.

Another example of Gabriel-type algebra is radical-graded artinian
algebra as mentioned in Theorem 2.4.

For a Gabriel-type algebra,  as Theorem 3.5 and 4.4 in \cite{Li},
 the uniqueness of the correspondent generalized path algebra and
 quiver holds up to isomorphism in the case the ideal is admissible, that is, if there exists another quiver and its
 related generalized path algebra such that the same isomorphism relation as in Definition 3.1 is satisfied for some admissible ideal $I$,
 then this quiver and the related
 generalized path algebra are just respectively the natural quiver and the corresponding one of
 this algebra. Exactly, we have the
following statement on the uniqueness:
\begin{theorem}
Assume $A$ is an artinian  algebra, $r=r(A)$ is the radical of
$A$. Let $A/r=\bigoplus^{p}_{i=1}\overline{A}_{i}$ with simple
ideals $\overline{A}_{i}$. If there is a quiver $Q$ and a normal
generalized path algebra $k(Q,\mathcal{B})$ with a set of simple
algebras $\mathcal{B}=\{B_{1},\cdot\cdot\cdot,B_{q}\}$ and an
admissible ideal $I$  of $k(Q,\mathcal{B})$ (i.e. for some $s$,
$J^{s}\subset I\subset J^2$)  such that $A\cong
k(Q,\mathcal{B})/I$ where $J$ the ideal of $k(Q,\mathcal{B})$
generated by all $\mathcal B$-paths of length one, then $Q$ is
just the natural quiver $\Delta_A$ of $A$ and $p=q$ such that
$\overline{A}_{i}\cong B_{i}$ for $i=1,...,p$ after reindexed. It
follows that $A$ is a Gabriel-type algebra.
\end{theorem}

{\em Proof:} Since $(J/I)^{s}\subseteq J^{s}/I=0$ and
$k(Q,\mathcal{B})/I/J/I\cong k(Q,\mathcal{B})/J=B_{1}\oplus \cdots
\oplus B_{q}$ semisimple, then rad$(k(Q,\mathcal{B})/I)=J/I$. From
the isomorphism $A\cong k(Q,\mathcal{B})/I$, we have
$A/$rad$A$$\cong k(Q,\mathcal{B})/I/J/I$, i.e. $\overline
A_{1}\oplus \cdots \oplus \overline A_{p}\cong B_{1}\oplus \cdots
\oplus B_{q}$. Thus $p=q$ and $\overline{A}_{i}\cong B_{i}$ for
$i=1,...,p$ after reindexed.

 By the isomorphism, the two
algebras $A$ and $k(Q,\mathcal{B})/I$ have the same natural
quivers, i.e. $\Delta_A=\Delta$. Then we only need to show that
the natural quiver $\Delta$ of $k(Q,\mathcal{B})/I$ is just $Q$.
Firstly since $p=q$, $\Delta_{0}=Q_{0}$. And the number of arrows
from $i$ to $j$ in $\Delta_{1}$ is $rank(B_{j}(J/I/J^{2}/I)B_{i})=
rank(B_{j}(J/J^{2})B_{i})$, which is just the number of arrows
from $i$ to $j$ in $Q_{1}$. Therefore $Q=\Delta=\Delta_{A}$.
$\;\;\;\square$

\begin{lemma}Let $A$ be a Gabriel-type artinian algebra with $A\stackrel\pi\cong k(\Delta_A,\mathcal A)/I$ for an ideal $I$ of $k(\Delta_A,\mathcal A)$
satisfying $I\subset J$. Assume that
$\{\varepsilon_i:\;i\in(\Delta_A)_0\}$ is
 the complete set of
non-isomorphic primitive orthogonal idempotents of $A$. Then,
there is a complete set of non-isomorphic primitive orthogonal
idempotents $\{d_i:\;i\in(\Delta_A)_0\}$ of $A/r$ such that
$\pi(\varepsilon_i)=d_i+I$ for any $i\in(\Delta_A)_0$.
\end{lemma}
{\em Proof}: Let $\pi(\varepsilon_i)=\widetilde\varepsilon_i+I$,
 then $\{\widetilde\varepsilon_i+I:\;i\in(\Delta_A)_0\}$ is a complete set of non-isomorphic primitive orthogonal idempotents of $k(\Delta_A,\mathcal A)/I$
 since $\pi$ is an isomorphism.

 Since $(\widetilde\varepsilon_i+I)^2=\widetilde\varepsilon_i+I$, we get
$\widetilde\varepsilon_i^2-\widetilde\varepsilon_i\in I$. Note
that $I$ lies in the ideal of $k(\Delta_A,\mathcal A)$ generated
by all $\mathcal A$-paths of length one. Because the square of any
non-cyclic path is zero, either
$\widetilde\varepsilon_i+I=E_ic_i+I$ or
$\widetilde\varepsilon_i+I=d_i+I$ where $E_i$ are  circles in
$\Delta_A$, $c_i$ and $d_i$ are primitive idempotents in
$k(\Delta_A,\mathcal A)$, or equivalently in $A/r$.

Let $w=1_{k(\Delta_A,\mathcal
A)/I}-\sum_{l\in(\Delta_A)_0}(\widetilde\varepsilon_i+I)$, then
$w$ is an idempotent and can be decomposed into a sum of some
primitive orthogonal idempotents $\widetilde f_j+I$, write
$w=(\widetilde f_1+I)+\cdots +(\widetilde f_t+I)$.  Thus,
$1_{k(\Delta_A,\mathcal A)/I}=\sum_{i\in(\Delta_A)_0}(\widetilde
\varepsilon_i+I)+\sum_{j=1}^t(\widetilde f_j+I).$

 Let $X+I$ and
$Y+I$ denote the sums of those idempotents in $\{\widetilde
\varepsilon_i+I:\;i\in(\Delta_A)_0\}\bigcup\{\widetilde
f_j+I:\;j=1+\cdots+t\}$ respectively in the forms $E_pc_p+I$ and
$d_q+I$, where $c_p, d_q\in A/r$. Thus,  $1_{k(\Delta_A,\mathcal
A)}+I=1_{k(\Delta_A,\mathcal A)/I}=(X+I)+(Y+I)$, it follows that
$X+I=(1_{k(\Delta_A,\mathcal A)}-Y)+I$.

 Suppose there are some $i$ such
that $\widetilde\varepsilon_i+I=E_ic_i+I\not=0$. Then $X+I\not=0$.
Hence $1_{k(\Delta_A,\mathcal A)}-Y\not=0$, then
$1_{k(\Delta_A,\mathcal A)}-Y\in X+I\subset J$, which is
impossible due to $1_{k(\Delta_A,\mathcal A)}-Y\in
k((\Delta_A)_0,\mathcal A)$.

The above contradiction means that each
$\widetilde\varepsilon_i+I=d_i+I$ where each $d_{i}$ is primitive
idempotent in $A/r$.

Clearly $\{d_i:\;i\in(\Delta_A)_0\}$ is a set of non-isomorphic
primitive orthogonal idempotents of $A/r$, by Lemma 3.1 it is a
complete set of non-isomorphic primitive orthogonal idempotents of
$A/r$. $\;\;\;\square$

\begin{theorem}
Let $A$ be a Gabriel-type artinian algebra over a field $k$  with
$A\stackrel\pi\cong k(\Delta_A,\mathcal A)/I$ for an ideal $I$ of
$k(\Delta_A,\mathcal A)$ satisfying $I\subset J$.
 Then for the basic algebra $B$ of $A$ and the basic algebra $C$ of $k(\Delta_A,\mathcal A)$, it holds that $B\cong (C+I)/I$.
\end{theorem}

{\em Proof:}$\;$ Let $\{\varepsilon_i:\;i\in(\Delta_A)_0\}$ be a
complete set of non-isomorphic primitive orthogonal idempotents of
$A$. Then, by Lemma 3.4, there is a complete set of non-isomorphic
primitive orthogonal idempotents $\{d_i: i\in(\Delta_A)_0\}$ of
$A/r$ such that $\pi(\varepsilon_i)=d_i+I$ for each
$i\in(\Delta_A)_0$.

By Lemma 3.1, $\{d_i: i\in(\Delta_A)_0\}$ is also a complete set
of non-isomorphic primitive orthogonal idempotents of
$k(\Delta_A,\mathcal A)$. Thus, by Proposition 3.2, we have
 $C\cong\bigoplus_{i,j\in(\Delta_A)_0}d_{i}k(\Delta_A,\mathcal A)d_{j}$.
Moreover, under the isomorphism $\pi$,
\begin{eqnarray*}B&\cong &\bigoplus_{i,j\in(\Delta_A)_0}\varepsilon_iA\varepsilon_j\\
&\cong &
\bigoplus_{i,j\in(\Delta_A)_0}(d_{i}+I)(k(\Delta_A,\mathcal
A)/I)(d_{j}+I)\\
&\cong &
\bigoplus_{i,j\in(\Delta_A)_0}(d_{i}k(\Delta_A,\mathcal A)d_{j}+I)/I\\
&\cong & (C+I)/I.\;\;\;\square\end{eqnarray*}

This theorem mentions the relation between the two basic algebras
$B$ and $C$ which are both decided by the same artinian algebra
$A$.

In general, for a Gabriel-type artinian $A$ whose the ideal $I$ is
 admissible (even only with $I\subset J$),
 the two natural quivers $\Delta_B$ and $\Delta_C$ of the associated  basic
 algebras $B$ of $A$ and $C$ of $k(\Delta_A,\mathcal A)$ are not
 equal.
In fact,  although $B\cong (C+I)/I$, rad$B\cong($rad$C+I)/I$, in
general $B/$rad$B\not\cong C/$rad$C$ and
$($rad$B)/($rad$B)^2\not\cong($rad$C)/($rad$C)^2$.

\begin{proposition} For a Gabriel-type artinian algebra $A$ with
$A\stackrel\pi\cong k(\Delta_A,\mathcal A)/I$, if $I$ is an
admissible ideal, the natural quivers of $A$ and
$k(\Delta_A,\mathcal A)$ are the same, i.e.
$\Delta_A=\Delta_{k(\Delta_A,\mathcal A)}$.
\end{proposition}
{\em Proof}:$\;\;$
 Since $I$ is admissible, there is a positive integer $s$ such that $J^s\subset I\subset J^2$.  rad$A \cong J/I$ as
proved in Theorem 3.3.  And $A\cong k(\Delta_A,\mathcal A)/I$,
then $A/$rad$A\cong k(\Delta_A,\mathcal A)/J$. Moreover,
rad$A/($rad$A)^2\cong J/J^2$. Thus, by the definition,
$\Delta_A=\Delta_{k(\Delta_A,\mathcal A)}$. $\;\;\;\square$

In the other case, for an artinian algebra $A$, when  $\Delta_A$
is admissible (i.e. is acyclic), it is true that
$\Delta_A=\Delta_{k(\Delta_A,\mathcal A)}$, since $J$ is just the
radical of $k(\Delta_A,\mathcal A)$.

To sum up, for a finite dimensional algebras $A$ over
algebraically closed field $k$,  when either  $\Delta_A$ is
admissible or $A$ is of Gabriel-type satisfying $A\cong
k(\Delta_A,\mathcal A)/I$ with admissible $I$, we have the
following diagram:

\[
\begin{array}{cclccclcc}
\Gamma_A&\supset &\Delta_A&=&\Delta_{k(\Delta_A,\mathcal A)}
&\subset&\Gamma_{k(\Delta_A,\mathcal A)}\\
\parallel &   &\cap &       &\cap &  &\parallel \\
\Gamma_B  & = & \Delta_B&  & \Delta_C & =& \Gamma_C
\end{array}
\]
where $\Gamma_A$ is the ordinary quiver of $A$, etc.; $\subset$,
$\supset$ and $\cap$ mean the embeddings of the dense
sub-quivers.

We feel the relations in this diagram would still hold for any
artinian algebras. This point of view will be discussed in the
subsequent work.

As we say above, the ordinary quiver and the natural quiver of a
finite dimensional  basic algebra coincide each other. In the end
of this section, we give an example which means the coincidence is
also possible to happen for some non-basic algebras. Meanwhile, in
this example, we show a method of computing the number of arrows
of the natural quiver of an artinian algebra.
\\
\\
 {\bf Example}$\;$ Let $k$ be an algebraically closed field of
 characteristic different from $2$ and let $Q$ be the quiver:
 \begin{figure}[hbt]
\begin{picture}(50,10)(-120,-10)
\put(100,0){\makebox(0,0){$ \bullet$}}
\put(102,-10){\makebox(0,0){$e_1$}}
 \put(105,0){\vector(1,0){30}}
\put(95,0){\vector(-1,0){30}}
\put(140,0){\makebox(0,0){$\bullet$}}
\put(144,-10){\makebox(0,0){$e_{2'}$}}
\put(145,0){\vector(1,0){30}}
\put(180,0){\makebox(0,0){$\bullet$}}\put(190,0){\makebox(0,0){$e_{3'}$}}
 \put(60,0){\makebox(0,0){$\bullet$}}
 \put(62,-10){\makebox(0,0){$e_2$}}
\put(55,0){\vector(-1,0){30}}\put(40,5){\makebox(0,0){$\beta$}}\put(85,5){\makebox(0,0){$\alpha$}}
\put(120,5){\makebox(0,0){$\alpha'$}}\put(160,5){\makebox(0,0){$\beta'$}}
\put(20,0){\makebox(0,0){$\bullet$}}\put(12,0){\makebox(0,0){$e_3$}}
\end{picture}
\end{figure}

Denote the path algebra $kQ$ by $\Lambda$ and let
$G=\langle\sigma\rangle$ be the group of order $2$. For the
elements $e_1, e_2, e_3, \alpha, \beta$ in $\Lambda$, let $\sigma
e_1=e_1, \sigma e_2=e_{2'}, \sigma e_3=e_{3'},
\sigma\alpha=\alpha', \sigma\beta=\beta'$. Then, there is only one
way of extending $\sigma$ to a $k$-algebra automorphism of $Q$ and
this is the way we will consider $G$ as a group of automorphisms
of $Q$. Now, we consider the ordinary quiver and the natural
quiver of the skew group algebra $\Lambda G$ (see \cite A).

Let $r$ be the radical of $\Lambda$. By Proposition 4.11 in \cite
A, $r\Lambda G=rad(\Lambda G)$. It is easy to see that $(\Lambda
G)/(r\Lambda G)\cong (\Lambda/r)G$. In the page 84 of \cite A, it
was given that $(\Lambda/r)G\cong A_1\times A_2\times A_3\times
A_4=k\times k\times\left(
\begin{array}{cl}
k & k \\
k & k
\end{array}\right)\times\left( \begin{array}{cl}
k & k \\
k & k
\end{array}\right)$ as algebras and the associated basic
algebra $B$ is obtained in the reduced form from $\Lambda G$,
which is Mortia-equivalent to $\Lambda G$, and moreover, it was
proved in \cite A that $B$ is isomorphic to the path algebra of
the following quiver:

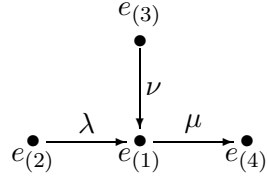
\begin{figure}[hbt]
\begin{picture}(50,50)(-120,-10)
\put(100,0){\makebox(0,0){$ \bullet$}}
\put(100,-7){\makebox(0,0){$e_{(1)}$}}\put(60,-7){\makebox(0,0){$e_{(2)}$}}
 \put(105,0){\vector(1,0){30}} \put(100,35){\vector(0,-1){30}}
\put(100,38){\makebox(0,0){$\bullet$}}\put(105,20){\makebox(0,0){$\nu$}}\put(100,48){\makebox(0,0){$e_{(3)}$}}
\put(65,0){\vector(1,0){30}} \put(140,0){\makebox(0,0){$\bullet$}}
\put(140,-7){\makebox(0,0){$e_{(4)}$}}
 \put(60,0){\makebox(0,0){$\bullet$}}\put(80,7){\makebox(0,0){$\lambda$}}\put(120,7){\makebox(0,0){$\mu$}}
\end{picture}
\end{figure}
.

This quiver is just the ordinary quiver $\Gamma_{\Lambda G}$ of
$\Lambda G$. Therefore, all $m_{ij}=0$, or $1$.

 For $i=1,2,3,4,\;dim_kA_i=n_i^2$ where $n_1=n_2=1$, $n_3=n_4=2$. By definition, for $i,j=1,2,3,4$, $t_{ij}$ is the
 rank of $_jM_i=A_j(r\Lambda G)/(r\Lambda G)^2A_i$ as
$A_j$-$A_i$-bimodule, equivalently, as a right
$\overline{A}_i\otimes\overline{A}_j^{op}$-module.  $A_i\otimes
A_j^{op}\cong M_{n_in_j}(k)$ is a simple algebra with dimension
$n_i^2n_j^2$. Thus, $_jM_i$ is semisimple over this simple
algebra. Let $_jM_i=L_1\oplus\cdots\oplus L_s$ where all $L_v$ are
simple $A_i\otimes A_j^{op}$-modules for $v=1,\cdots,s$. $L_v$ can
be considered as a simple right ideal of $A_i\otimes A_j^{op}$,
therefore, $L_v\cong(k\;k\;\cdots\;k)$ the whole set of all
$1\times n_in_j$ matrices over $k$. For any $0\not =x_v\in L_v$,
$L_v=x_v(A_i\otimes A_j^{op})$. Then,
$_jM_i=\bigoplus_{v=1}^sx_v(A_i\otimes A_j^{op})$
 as $(A_i\otimes A_j^{op})$-modules.

First, we prove $s=m_{ij}$ the number of the arrows from  the
vertex $i$ to the other vertex $j$ in the ordinary quiver
$\Gamma_{\Lambda G}$ of $\Lambda G$.

Let  $\{S_{1},S_{2},\cdots, S_{n}\}$ be  the complete set of
non-isomorphic simple $\Lambda G$-modules. Then
$m_{ij}=dim_kExt_A(S_i,S_j)$. By \cite{Liu}\cite{A},
$dim_kExt_A(S_i,S_j)=dim_k(k\varepsilon_j(r\Lambda G)/(r\Lambda
G)^2k\varepsilon_i)$ where $\{\varepsilon_{1},\cdots,
\varepsilon_{n}\}$ is the complete set of primitive orthogonal
idempotents of $(\Lambda G)/rad(\Lambda G)$ with $\varepsilon_i\in
 A_i$. And,
  $k\varepsilon_j(r\Lambda
G)/(r\Lambda G)^2k\varepsilon_i= \varepsilon_jA_j(r\Lambda
G)/(r\Lambda G)^2A_i\varepsilon_i
=\varepsilon_j$$_jM_i\varepsilon_i\cong\;_jM_i(\varepsilon_i\otimes
\varepsilon_j)=(L_1\oplus\cdots\oplus L_s)(e_i\otimes
e_j)=(x_1\oplus\cdots\oplus
 x_s)(A_i\otimes A_j^{op})(e_i\otimes e_j)\cong(x_1\oplus\cdots\oplus
 x_s)M_{n_in_j}(k)E_{ll}$ where $\varepsilon_i\otimes\varepsilon_j$ is a primitive
 idempotent of $A_i\otimes A_j^{op}$  so let $E_{ll}$ be the correspondent element of $\varepsilon_i\otimes
 \varepsilon_j$ in $M_{n_in_j}(k)$ under the isomorphism.

  Obviously,
 $dim_kx_vM_{n_in_j}(k)E_{ll}=1$ for all $v=1,\cdots,s$. Then, $dim_k(x_1\oplus\cdots\oplus
 x_s)M_{n_in_j}(k)E_{ll}=s$. It follows that
  $m_{ij}=s$.

For each pair $(i,j)$, when $m_{ij}=0$, we have
$k\varepsilon_j(r\Lambda G)/(r\Lambda G)^2k\varepsilon_i=0$. Then
$_jM_i=0$. Thus, the rank $t_{ij}$ of $_jM_i$ equals $0$. When
$m_{ij}=1$, then $s=m_{ij}=1$, that is, $_jM_i=L_1$ is a simple
$(A_i\otimes A_j^{op})$-module. Hence, the rank $t_{ij}$ of
$_jM_i$ is 1 in this case.

According to the above discussion, for each pair $(i,j)$, we have
$t_{ij}=m_{ij}=0$ or 1. Therefore, the natural quiver
$\Delta_{\Lambda G}$ is equal to the ordinary quiver
$\Gamma_{\Lambda G}$.

\section{Interpretations}

In \cite{A}\cite{ASS}, given a finite dimensional algebra $A$, the
ordinary quiver $\Gamma_A$ can be constructed
 by the indecomposable projective modules
and the irreducible morphisms between them.
 So the ordinary quiver of $A$ provides a
convenient way to study its projective (or injective) modules and
morphisms between them, even when $A$ is not a basic algebra.  By
the Gabriel theorem,  the ordinary quiver of a finite-dimensional
algebra $A$ is used as a tool to characterize the structure of its
associated basic algebra but not of $A$. In this reason, the
ordinary quiver is not effective enough to characterize a
non-basic algebra.  The generalized Gabriel theorem in \cite{Li}
shows the arrival of our goal via the natural quiver under some
conditions.

Note that the AR-quiver of the sub-category of {\bf proj}$A$ with
irreducible morphisms is isomorphic to the opposite of the
ordinary quiver of $A$.

Through \cite{Li} and here, we think the method of natural quiver
may offset some shortage of ordinary quiver and AR-quiver. In
certain sense,  the natural quiver of an artinian algebra $A$ will
also be available for the theory of representations of an artinian
algebra.

Under certain condition, the representation category {\bf Rep}$A$
of $A$ can be decided wholly by the ordinary quiver and the
AR-quiver.
 The category of representations of $A$ may be partially induced from the category of representations of $\Gamma_B$ through the basic algebra
 $B$. For example,
  when $A$ is Gabriel-type, that is,  $A$ is isomorphic to some
  quotient of the generalized path algebra of $\Delta_A=\Delta_B$,
  any representations of $A$ can be induced directly from some of representations of the
  generalized path algebra of $\Delta_A$. In the classical theory of representations of artin algebras (see \cite A\cite{ASS}\cite{D} etc.), one wants to characterize
   {\bf Rep}$A$ through representations of $\Delta_B$ with $B$. However, the difficulty is that in general, it is not easy to
   construct concretely the  basic algebra $B$ from $A$.
 By comparison,  the method of natural quivers is more straightforward through representations of the
  generalized path algebra of $\Delta_A$. Therefore, we hope to set up this new approach to representations of an artinian algebra via representations of the
  generalized path algebra of its natural quiver.

\end{document}